# On $g$-functions for subshifts

Wolfgang Krieger[1]

*University of Heidelberg*

**Abstract:** A necessary and sufficient condition is given for a subshift presentation to have a continuous $g$-function. An invariant necessary and sufficient condition is formulated for a subshift to posses a presentation that has a continuous $g$-function.

## 1. Introduction

Let $\Sigma$ be a finite alphabet, and let $S$ denote the shift on $\Sigma^{\mathbb{Z}}$,

$$S((x_i)_{i\in\mathbb{Z}}) = ((x_{i+1})_{i\in\mathbb{Z}}), \qquad (x_i)_{i\in\mathbb{Z}} \in \mathbb{Z}.$$

A closed $S$-invariant set $X \subset \Sigma^{\mathbb{Z}}$ with the restriction of $S$ acting on it, is called a subshift. A finite word is said to be admissible for a subshift if it appears in a point of the subshift. A subshift is uniquely determined by its set of admissible words. A subshift is said to be of finite type if its admissible words are defined by excluding finitely many words from appearing as subwords in them. Subshifts are studied in symbolic dynamics. For an introduction to symbolic dynamics see [7] and [10].

We introduce notation. Given a subshift $X \subset \Sigma^{\mathbb{Z}}$ we set

$$x_{[i,k]} = (x_j)_{i\leq j\leq k}, \quad x \in X, i,k \in \mathbb{Z}, i \leq k,$$

and

$$X_{[i,k]} = \{x_{[i,k]} : x \in X\}.$$

We use similar notation also for blocks,

$$b_{[i',k']} = (b_j)_{i'\leq j\leq k'}, \quad b \in X_{[i,k]},\ i \leq i' \leq k' \leq k,$$

and also if indices range in semi-infinite intervals. Blocks also stand for the words they carry. We denote

$$\Gamma_n^+(x^-) = \{b \in X_{[1,n]} : (x^-, b) \in X_{(-\infty,n]}\}, \qquad n \in \mathbb{N},$$

$$\Gamma^+(x^-) = \bigcup_{n\in\mathbb{N}} \Gamma_n^+(x^-),$$

$$\Gamma_\infty^+(x^-) = \{x^+ \in X_{[1,\infty)} : (x^-, x^+) \in X\}, \qquad x^- \in X_{(-\infty,0]}$$

and

$$\Gamma_n^+(a) = \{b \in X_{[1,n]} : (a,b) \in X_{(-k,n]}\}, \qquad n \in \mathbb{N},$$

$$\Gamma^+(a) = \bigcup_{n\in\mathbb{N}} \Gamma_n^+(a),$$

---

[1]Institute for Applied Mathematics, University of Heidelberg, Im Neuenheimer Feld 294, 69120 Heidelberg, Germany, e-mail: krieger@math.uni-heidelberg.de

$$\Gamma_\infty^+(a) = \{x^+ \in X_{[1,\infty)} : (a, x^+) \in X_{(-k,\infty)}\}, \qquad a \in X_{(-k,0]}, k \in \mathbb{Z}_+,$$

$\Gamma^-$ has the time symmetric meaning. We denote

$$\omega_n^+(a) = \bigcap_{x^- \in \Gamma_\infty^-(a)} \{b \in X_{[1,n]} : (x^-, a, b) \in X_{(-\infty,n]}\},$$

$$\omega^+(a) = \bigcup_{n \in \mathbb{N}} \omega_n^+(a), \qquad a \in X_{(-k,0]}, k \in \mathbb{Z}_+.$$

The notions of $g$-function and $g$-measure go back to Mike Keane's papers [5], [6]. Subsequently a substantial theory of $g$-functions and $g$-measures developed with contributions from many sides (see e.g. [1],[4],[14],[16],[17] and the references given there. For the origin of these notions see also [2]). These notions have formulations for general subshifts (see [11, p. 24]). We are interested in continuous $g$-functions and therefore introduce a $g$-function for a subshift $X \subset \Sigma^\mathbb{Z}$ as a continuous mapping

$$g : \{(x^-, \sigma) \in X_{(-\infty,0]} \times \Sigma : \sigma \in \Gamma_1^+(x^-)\} \to [0,1]$$

such that

$$\sum_{\alpha \in \Gamma_1^+(x^-)} g(x^-, \alpha) = 1, \qquad x^- \in X_{(-\infty,0]},$$

and a $g$-measure as an invariant probability measure $\mu$ of the subshift $X$ such that

$$\mu(\{x \in X : x_{[-k,1]} = (a, \alpha)\})$$
$$= \int_{\{x \in X : x_{[-k,0]} = a\}} g(x^-, \alpha) d\mu, \quad (a, \alpha) \in X_{[-k,1]}, k \in \mathbb{N}.$$

(Note that we have reversed the time direction.) We show in Section 2 that a subshift that has a strictly positive $g$-function is of finite type. Denote for $x^- \in X_{(-\infty,0]}$

$$\Delta_1^+(x^-) = \bigcup_{n \in \mathbb{N}} \omega_1^+(x^-_{[-n,0]}).$$

In Section 2 we prove that a subshift $X \subset \Sigma^\mathbb{Z}$ has a $g$-function if and only if for all $x^- \in X_{(-\infty,0]}$, $\Delta_1^+(x^-) \neq \emptyset$. We refer to this property of a subshift presentation as property $g$.

A directed graph with vertex set $\mathcal{M}$ and edges carrying labels taken from a finite alphabet $\Sigma$ is called a Shannon graph if the labeling is 1-right resolving in the sense that for all $\mu \in \mathcal{M}$ and $\sigma \in \Sigma$ there is at most one edge leaving $\mu$ that carries the label $\sigma$. Denote here the set of initial vertices of the edges that carry the label $\sigma$ by $\mathcal{M}(\sigma)$, and for $\mu \in \mathcal{M}(\sigma)$ denote by $\tau_\sigma(\mu)$ the final vertex of the edge that leaves $\mu$ and carries the label $\sigma$. The Shannon graph $\mathcal{M}$ is determined by the transition rules $(\tau_\sigma)_{\sigma \in \Sigma}$. A Shannon graph is said to present a subshift $X \subset \Sigma^\mathbb{Z}$ if every vertex has an edge leaving it and an edge entering it, and if the set of admissible words of the subshift coincides with the set of label sequences of finite paths in the graph.

For a finite alphabet $\Sigma$ denote by $\mathcal{M}(\Sigma)$ the set of probability measures on $\Sigma^\mathbb{N}$ with its weak *-topology. With the notation

$$C(a) = \bigcap_{1 \leq i \leq n} \{(x_i)_{i \in \mathbb{N}} \in \Sigma^\mathbb{N} : a_i = x_i\}, \qquad (a_i)_{1 \leq i \leq n} \in \Sigma^N, n \in \mathbb{N},$$

$$\mathcal{M}(\Sigma)(\sigma) = \{\mu \in \mathcal{M}(\Sigma) : \mu(C(\sigma)) > 0\}, \qquad \sigma \in \Sigma,$$

let for $\mu \in \mathcal{M}(\Sigma)(\sigma), \tau_\sigma(\mu)$ be equal to the conditional measure of $\mu$ given $C(\sigma)$,

$$\tau_\sigma(\mu)(C(b)) = \frac{\mu(C(\sigma,b))}{\mu(C(\sigma))}, \qquad b \in \Sigma^N, N \in \mathbb{N}.$$

In this way $\mathcal{M}(\Sigma)$ has been turned into a Shannon graph with the transition rules $(\tau_\sigma)_{\sigma \in \Sigma}$. The Shannon graph $\mathcal{M}(\Sigma)$ is accompanied by another Shannon graph with vertex set $\bigcup_{N \in \mathbb{Z}_+} \mathcal{M}_N(\Sigma)$, where for $N \in \mathbb{N}$, $\mathcal{M}_N(\Sigma)$ is the set of probability vectors on $\Sigma^N$, and where $\mathcal{M}_0 = \{\emptyset\}$. With the notation

$$C_N(\sigma) = \{(a_i)_{1 \leq i \leq n} \in \Sigma^N : a_1 = \sigma\},$$

$$\mathcal{M}_N(\Sigma)(\sigma) = \{\mu \in \mathcal{M}_N(\Sigma) : \mu(C_N(\sigma)) > 0\},$$

one sets for $\sigma \in \Sigma, \mu \in \mathcal{M}_N(\Sigma)(\sigma), N > 1, \tau_\sigma(\mu)$ equal to the probability vector $\nu \in \mathcal{M}_{N-1}(\Sigma)$ that is given by

$$\nu(b) = \frac{\mu(\sigma,b)}{\mu(C(\sigma))}, \qquad b \in \Sigma^{N-1},$$

and one sets $\tau_\sigma(\mu) = \{\emptyset\}$ for $\sigma \in \Sigma, \mu \in \mathcal{M}_1(\Sigma)$. In this way $\bigcup_{N \in \mathbb{Z}_+} \mathcal{M}_N$ has been turned into a Shannon graph that one can equip further with the restriction mapping $\iota$ that assigns to $\mu \in \mathcal{M}_N(\Sigma), N > 1$, its marginal vector in $\mathcal{M}_{N-1}(\Sigma)$, and that assigns to a $\mu \in \mathcal{M}_1(\Sigma)$ the empty set. The mapping $\iota$ commutes with the transition rules of the Shannon graph. Call a set $\mathcal{M} \subset \mathcal{M}(\sigma)$ transition complete if for $\sigma \in \Sigma, \mu \in \mathcal{M} \cap \mathcal{M}(\Sigma)$ implies that also $\tau_\sigma(\mu) \in \mathcal{M}$. Call a set $\mathcal{M} \subset \mathcal{M}(\sigma)$ in-complete if for all $\mu \in \mathcal{M}$ there is a $\nu \in \mathcal{M}$ that is the initial vertex of an edge ending in $\mu$. Every transition complete and in-complete set $\mathcal{M} \subset \mathcal{M}(\sigma)$ determines a Shannon graph with transition rules that are inherited from the Shannon graph $\mathcal{M}(\Sigma)$. These sub-Shannon graphs $\mathcal{M} \subset \mathcal{M}(\sigma)$ are accompanied by sub-Shannon graphs $\bigcup_{N \in \mathbb{Z}_+} \mathcal{M}_N \subset \bigcup_{N \in \mathbb{Z}_+} \mathcal{M}_N(\Sigma)$ where $\mathcal{M}_N$ contains the probability vectors that are given by the marginals of the measures in $\mathcal{M}$, and where the transition rules and the mapping $\iota$ are passed down from $\bigcup_{N \in \mathbb{Z}_+} \mathcal{M}_N(\Sigma)$. In [12] Kengo Matsumoto introduced a class of structures that he called $\lambda$-graph systems. $\lambda$-graph systems have the form of a Bratteli diagram, that is, they have a finite number of vertices at each level. In the structures $\bigcup_{N \in \mathbb{Z}_+} \mathcal{M}_N$ the sets $\mathcal{M}_N, N \in \mathbb{N}$, are not necessarily finite, but otherwise these structures have all the attributes of a $\lambda$-graph system. We will refer to them as measure $\lambda$-graph systems. We say that a measure $\lambda$-graph system presents a subshift if the set of admissible words of the subshift coincides with the set of label sequences of finite paths in the measure $\lambda$-graph system.

In Section 3 we are concerned with the measure $\lambda$-graph system that is generated by a $g$-function $g$ of a subshift $X \subset \Sigma^\mathbb{Z}$. The continuity of the $g$-function translates into a property of the generated measure $\lambda$-graph system that we call contractivity. Every contractive measure $\lambda$-graph system determines a $g$-function of the subshift that it presents and it is in turn generated by this $g$-function. A subshift that is presented by a contractive measure $\lambda$-graph system has a property that we call property $(D)$. In Section 4 we prove the invariance of property $(D)$ under topological conjugacy, and point to some classes of subshifts that have property $(D)$ and that have presentations with property $g$. Every subshift that has property $(D)$ and that has a presentation with property $g$ admits a presentation by a contractive measure $\lambda$-graph system.

## 2. g-functions of subshifts

**Lemma 2.1.** *Let $X \subset \Sigma^{\mathbb{Z}}$ be a subshift, let $Y^- \subset X_{(-\infty,0]}$ be dense in $X_{(-\infty,0]}$, and let*
$$g : \{(x^-, \sigma) \in Y^- \times \Sigma : \sigma \in \Gamma_1^+(x^-)\} \to [0,1]$$
*be a continuous mapping such that*
$$g(y^-, \alpha) = 0, \qquad y^- \in Y^-, \alpha \notin \Gamma_1^+(y^-),$$
*and such that*
$$\sum_{\beta \in \Gamma_1^+(y^-)} g(y^-, \beta) = 1, \qquad y^- \in Y^-. \tag{1}$$
*Let*
$$x^- \in Y^-, \quad \alpha \in \Gamma_1^+(x^-) \setminus \Delta_1^+(x^-). \tag{2}$$
*Then*
$$g(x^-, \alpha) = 0. \tag{3}$$

*Proof.* By (2) there are $M_n \in \mathbb{N}$, $n \in \mathbb{N}$, and $a(n) \in X_{[-M_n,0]}$, such that
$$a(n)_{[-n,0]} = x_{[-n,0]} \tag{4}$$
and
$$\alpha \notin \Gamma_1^+(a(n)), \qquad n \in \mathbb{N}. \tag{5}$$

Since $Y^-$ is dense in $X_{(-\infty,0]}$ one can find $y^-(n) \in Y^-$ such that
$$y^-(n)_{[-M_n,0]} = a(n), \qquad n \in \mathbb{N}. \tag{6}$$

By (4) and (6)
$$\lim_{n \to \infty} y^-(n) = x^-, \tag{7}$$
and from (5) and (6)
$$\alpha \notin \Gamma_1^+(y^-(n)), \qquad n \in \mathbb{N}. \tag{8}$$

Choose an increasing sequence $n_k, k \in \mathbb{N}$, and a set $\Gamma \subset \Sigma$ such that
$$\Gamma_1^+(y^-(n_k)) = \Gamma, \qquad k \in \mathbb{N}. \tag{9}$$

Then by (7) and by compactness of $X$
$$\Gamma \subset \Gamma_1^+(x^-), \tag{10}$$
and from (1) and by the continuity of $g$
$$\sum_{\gamma \in \Gamma} g(x^-, \gamma) = 1,$$
and then (3) follows from (8) and (9). □

**Corollary 2.2.** *Let the subshift $X \subset \Sigma^{\mathbb{Z}}$ have a strictly positive g-function. Then $X$ is of finite type.*

*Proof.* For $\alpha \in \Sigma$ one has the closed set

$$X_\alpha^- = \{x^- \in X_{(-\infty,0]} : \alpha \in \Gamma_1^+(x^-)\}.$$

It follows from Lemma 2.1 that

$$\Gamma_1^+(x^-) = \Delta_1^+(x^-), \qquad x^- \in X_{(-\infty,0]}.$$

Therefore the increasing sequence

$$X_\alpha^-(n) = \{x^- \in X_{(-\infty,0]} : \alpha \in \omega_1^+(x^-_{[-n,0]})\}, n \in \mathbb{N},$$

of open subsets of $X_{(-\infty,0]}$ is a cover of $X_\alpha^-$. One has therefore an $n_\alpha \in \mathbb{N}$ such that

$$X_\alpha^- = X_\alpha^-(n_\alpha).$$

With

$$N = \max\{n_\alpha : \alpha \in \Sigma\},$$

one has

$$\Gamma_1^+(x^-) \subset \omega_1^+(x^-_{[-N,0]}), \qquad x^- \in X_{(-\infty,0]},$$

which means that the subshift $X$ is determined by its set of admissible words of length $N+2$. $\square$

**Lemma 2.3.** *Let $X \subset \Sigma^\mathbb{Z}$ be a subshift such that*

$$\{x^- \in X_{(-\infty,0]} : \Delta_1^+(x^-) \neq \emptyset\} \neq \emptyset.$$

*Then there exists a continuous mapping*

$$g : \{x^- \in X_{(-\infty,0]} : \Delta_1^+(x^-) \neq \emptyset\} \times \Sigma \to [0,1]$$

*such that*

$$g(x^-, \alpha) = 0, \qquad \alpha \notin \Delta_1^+(x^-),$$

*and*

$$g(x^-, \alpha) > 0, \qquad \alpha \in \Delta_1^+(x^-).$$

*Proof.* For $x^- \in X_{(-\infty,0]}$ and $\alpha \in \Delta_1^+(x^-)$ set

$$n(x^-, \alpha) = \min\{n \in \mathbb{N} : \alpha \in \omega^+(x^-_{[-n,0]})\},$$

and

$$n(x^-) = \min_{\alpha \in \Delta_1^+(x^-)} n(x^-, \alpha).$$

For $x^- \in X_{(-\infty,0)}$, and $\gamma \in \Gamma_1^+(x^-)$, set

$$g(x^-, \gamma) = \begin{cases} 0, & \text{if } \gamma \notin \Delta_1^+(x^-), \\ \dfrac{(n(x^-,\gamma)-n(x^-))^{-1}}{\sum_{\beta \in \Delta^+(x^-)}((n(x^-,\beta)-n(x^-))^{-1}}, & \text{if } \gamma \in \Delta_1^+(x^-). \end{cases}$$

To prove continuity of the mapping $g$ at a point $(x^-, \alpha)$, , $x^- \in X_{(-\infty,0]}, \alpha \in \Delta^+(x^-)$, let

$$N(x^-) = \max_{\beta \in \Delta_1^+(x^-)} n(x^-, \beta),$$

and let $y^-(k) \in X_{(-\infty,0]}, k \in \mathbb{N}$, be such that

$$\lim_{k \to \infty} y^-(k) = x^-.$$

For $M \in \mathbb{N}$, let $k_\circ \in \mathbb{N}$ be such that

$$y^-(k)_{[-M-N[x^-),0]} = x^-_{[-M-N[x^-),0]}, \qquad k \geq k_\circ.$$

Then

$$\alpha \in \Delta^+(y^-(k)), \qquad k \geq k_\circ,$$

and

$$|g(y^-(k), \alpha) - g(x^-, \alpha)| < M^{-1} |\Sigma| g(x^-, \alpha), \qquad k \geq k_\circ.$$

To prove continuity of the mapping $g$ at a point $(x^-, \alpha), x^- \in X_{(-\infty,0]}, \Delta_1^+(x^-) \neq \emptyset, \alpha \notin \Delta_1^+(x^-)$, let $y^-(k) \in X, k \in \mathbb{N}$, be such that

$$\lim_{k \to \infty} y^-(k) = x^-,$$

and such that

$$\alpha \in \Delta_1^+(y^-(k)).$$

For $M \in \mathbb{N}$, let $k_\circ \in \mathbb{N}$ be such that

$$y^-(k)_{[-M-n(x^-),0]} = x^-_{[-M-n(x^-),0]}, \qquad k \geq k_\circ.$$

Then

$$g(y^-(k), \alpha) < \frac{1}{M}, \qquad k \geq k_\circ. \qquad \square$$

A $g$-function of a subshift such that $g(x^-, \alpha) > 0$ for $x^- \in X_{(-\infty,0]}, \alpha \in \Delta_1^+(x^-)$ we will call a strict $g$-function.

**Theorem 2.4.** *The following are equivalent for a subshift $X \subset \Sigma^{\mathbb{Z}}$: (a) $X$ has a $g$-function. (b) $X$ has property g. (c) $X$ has a strict $g$-function.*

*Proof.* That (a) implies (b) follows from Lemma 2.1. That (b) implies (c) follows from Lemma 2.3. $\square$

## 3. Presentations of subshifts and property (*D*)

Given a $g$-function of the subshift $X \subset \Sigma^{\mathbb{Z}}$ we define inductively probability vectors $\mu_n(x^-) \in \mathcal{M}_n(\Sigma), n \in \mathbb{N}, x^- \in X_{-\infty,0]}$, by setting $\mu_n(x^-)(a)$ equal to zero, if $a \in \Sigma^N$ is not in $\Gamma_N^+(x^-)$, and by setting

$$\mu_N(a) = \prod_{1 \leq k \leq N} g(x^-, a_{[1,k)}, a_k), \qquad a \in \Gamma_N^+(x^-), \ N \in \mathbb{N}, \tag{11}$$

and we let $\mu(x^-) \in \mathcal{M}(\Sigma)$ be the probability measure that has as marginal measures those that are given by the probability vectors $\mu_N(x^-), N \in \mathbb{N}$. We set

$$\mathcal{M}_N(X, g) = \{\mu_N(x^-) : x^- \in X_{(-\infty,0]}\}, \qquad N \in \mathbb{N},$$

$$\mathcal{M}(X, g) = \{\mu(x^-) : x^- \in X_{(-\infty,0]}\}.$$

Here $\bigcup_{N \in \mathbb{N}} \mathcal{M}_N(X, g) \subset \bigcup_{N \in \mathbb{N}} \mathcal{M}_N$ is the measure $\lambda$-graph system that accompanies the compact transition complete and in-complete sub-Shannon graph $\mathcal{M}(X, g)$ of $\mathcal{M}(\Sigma)$. Entities like the mapping that assigns to a point $x^- \in X_{(-\infty,0]}$ for a subshift $X \subset \Sigma^{\mathbb{Z}}$ the measure $\mu(x^-) \in \mathcal{M}(\Sigma)$, or the inverse image under this mapping of a single measure, appear prominently within a theory that was put forward by James Crutchfield et al (see e.g.[15]).

We set for a given transition complete and in-complete Shannon graph $\mathcal{M} \subset \mathcal{M}(\Sigma)$ and for its accompanying measure $\lambda$-graph system $\lambda$-graph system $\bigcup_{N \in \mathbb{Z}_+} \mathcal{M}_N$ inductively

$$\tau_a(\mu) = \tau_{a_{-m}}(\tau_{a_{(-m,0]}}(\mu)), \qquad 0 < m \leq n, a \in X_{[-n,0]}, n \in \mathbb{N}.$$

Call a compact transition complete Shannon graph $\mathcal{M} \subset \mathcal{M}(\Sigma)$ contractive if, with $X \subset \Sigma^{\mathbb{Z}}$ the subshift that is presented by $\mathcal{M}$, one has for all $x^- \in X_{(-\infty,0]}$ that the limits

$$\lim_{k \to \infty} \tau_{x^-_{[-n,0)}}(\mu), \qquad \mu \in \mathcal{M},$$

exist. Call a measure $\lambda$-graph system $\bigcup_{N \in \mathbb{Z}_+} \mathcal{M}_N$ contractive if, with $X \subset \Sigma^{\mathbb{Z}}$ the subshift that it presents one has that for all $x^- \in X_{(-\infty,0]}$, that

$$\lim_{k \to \infty} \operatorname{diam}( \bigcup_{\mu \in \mathcal{M}_n} \tau_{x^-_{[-n,0)}}(\mu)) = 0.$$

To a contractive Shannon graph there corresponds a contractive measure $\lambda$-graph system and vice versa. Due to this one-to-one correspondence between contractive measure $\lambda$-graph systems and contractive Shannon graphs one can formulate here arguments and results in terms of one or the other. We will express ourselves in terms of the contractive measure $\lambda$-graph systems, the motivation being, that the constituent elements of the contractive measure $\lambda$-graph systems are sequentially generated by the $g$-function according to (11). Note that the theory of contractive measure $\lambda$-graph systems and Matsumoto's theory of $\lambda$-graph systems intersect in the theory of topological Markov shifts.

**Proposition 3.1.** *Let $\bigcup_{n \in \mathbb{Z}_+} \mathcal{M}_N$ be a contractive measure $\lambda$-graph system. The subshift that is presented by $\bigcup_{n \in \mathbb{Z}_+} \mathcal{M}_N$ has a g-function $g$ such that $\mathcal{M}_N = \mathcal{M}_N(X, g), N \in \mathbb{Z}_+$.*

*Proof.* Every contractive measure $\lambda$-graph system $\bigcup_{n \in \mathbb{Z}_+} \mathcal{M}_n$ defines a $g$-function $g$ of the subshift $X \subset \Sigma^{\mathbb{Z}}$ that it presents by

$$(g(x^-, \alpha))_{\sigma \in \Gamma_1^+(x^-)} \in \bigcap_{n \in \mathbb{N}} \bigcup_{\mu \in \mathcal{M}_n} \tau_{x^-_{[-n,0)}}(\mu), \qquad x^- \in X_{(-\infty,0]}. \tag{12}$$

One uses the hypothesis that every vertex in $\bigcup_{n \in \mathbb{Z}_+} \mathcal{M}_N$ has a predecessor to show that the g-function that is associated to $\bigcup_{n \in \mathbb{Z}_+} \mathcal{M}_N$ according to (12) has the stated property. $\square$

We say that a subshift $X \subset \Sigma^{\mathbb{Z}}$ has property $(D)$ if for all admissible words $b\sigma$ of $X$ there exists a word $a \in \Gamma^-(b)$ of $X$ such that $\sigma \in \omega_1^+(ab)$.

**Lemma 3.2.** *A subshift $X \subset \Sigma^{\mathbb{Z}}$ that admits a presentation by a contractive $\lambda$-graph system $\bigcup_{n \in \mathbb{Z}_+} \mathcal{M}_N$ has property $(D)$.*

*Proof.* Let $b$ be an admissible word of the subshift $X$ and let $\sigma \in \Gamma_1^+(b)$. Denote the length of $b$ by $K$, and let $\mu \in \mathcal{M}_{K+1}$ be a vertex with a path leaving it that has label sequence $b$. Set $\nu = \tau_b(\mu)$. Then $\nu(\sigma) > 0$. Let then $x^- \in X_{(-\infty,0)}$ be such that $x^-_{[-K,0]} = b\sigma$, and such that $x^-_{(-\infty,K)}$ is the label sequence of a path in $\bigcup_{n \in \mathbb{Z}_+} \mathcal{M}_N$ that leads into the vertex $\mu$. Then it follows for the $g$-function that is associated to $\bigcup_{n \in \mathbb{Z}_+} \mathcal{M}_N$ according to (12) that $g(x^-, \sigma) = \nu(\sigma) > 0$. Apply Lemma 2.1 to conclude the proof. □

**Lemma 3.3.** *The following are equivalent for a subshift $X \subset \Sigma^{\mathbb{Z}}$: (a) $X$ has property $(D)$. (b) For all admissible words $a$ and $c$ of $X$ such that $c \in \Gamma^+(a)$ there exists an admissible word $b$ of $X$ such that $ac \in \omega^+(b)$.*

*Proof.* We prove that (a) implies (b). For this let $a$ and $c = (c_l)_{1 \leq l \leq k}, k \in \mathbb{N}$, be admissible words of the subshift $X$ such that $c \in \Gamma^+(a)$. Choose inductively words $b^l, 1 \leq l \leq k$, such that
$$c_l \in \omega_1^+((b^m)_{l \geq m \geq 1}, ac_{[1,l)}), \qquad 1 \leq l \leq k.$$
Then set
$$b = (b^l)_{k \geq l \geq 1}.$$
□

**Lemma 3.4.** *Let $X \subset \Sigma^{\mathbb{Z}}$ be a subshift with properties $g$ and $(D)$, and let $g$ be a strict $g$-function of $X$. Then the measure $\lambda$-graph system $\bigcup_{n \in \mathbb{Z}_+} \mathcal{M}_n(X, g)$ presents $X$.*

*Proof.* Let $a$ and $b$ be admissible words of the subshift $X$ such that $a \in \omega^+(b)$. Let $N$ denote the length of $a$, and let $K$ denote the length of $b$. Since the $g$-function $g$ is assumed strict it follows that for $x^- \in X_{(-\infty,0]}$ such that $x^-_{[-K,0]} = b$ one has that $\mu_N(x^-)(a) > 0$. This implies that there is a path in $\bigcup_{n \in \mathbb{Z}_+} \mathcal{M}_n(X, f)$ with label sequence $a$ that leads into the vertex $\emptyset$. □

**Theorem 3.5.** *The presentation of a subshift by a contractive measure $\lambda$-graph system has property g. A subshift that admits a presentation by a contractive measure $\lambda$-graph system has property $(D)$.*

*Proof.* The presentation of a subshift by a contractive measure $\lambda$-graph system has property g by Proposition 3.1 and a subshift that admits such a presentation has property $(D)$ by Lemma 3.2. Conversely, if a subshift presentation has property $g$, then by Lemma 2.3 it has a strict $g$-function and by Lemma 3.4 this strict $g$-function generates a contractive measure $\lambda$-graph system that presents the subshift. □

## 4. Invariance

We recall that, given subshifts $X \subset \Sigma^{\mathbb{Z}}, \bar{X} \subset \bar{\Sigma}^{\mathbb{Z}}$ and a continuous shift-commuting map $\varphi : X \to \bar{X}$ there is for some $L \in \mathbb{Z}_+$ a block mapping
$$\Phi : X_{[-L,L]} \to \bar{\Sigma}$$
such that
$$\varphi(x) = (\Phi(x_{[i-L,i+L]}))_{i \in \mathbb{Z}}.$$
We say then that $\varphi$ is implemented by $\Phi$, and we write
$$\Phi(a) = (\Phi(a_{[j-L,j+L]})_{i+L \leq j \leq k-L}), \quad a \in X_{[i,k]}, \quad k - i \geq 2L,$$
and use similar notation if indices range in semi-infinite intervals. Recall that the $n$-block system of a subshift $X \in \Sigma^{\mathbb{Z}}$ is its image in $(X_{[1,n]})^{\mathbb{Z}}$ under the mapping $x \to (x_{(i,i+n]})_{i \in \mathbb{Z}}, x \in X$.

Call a subshift $X \subset \Sigma^{\mathbb{Z}}$ right instantaneous [8] if for all $\sigma \in \Sigma, \omega_1^+(\sigma) \neq \emptyset$.

**Proposition 4.1.** *A subshift $X \subset \Sigma^{\mathbb{Z}}$ has a g-function if and only if one of its n-block systems is right-instantaneous.*

*Proof.* For $n \in \mathbb{N}$ and $\sigma \in \Sigma$ one has the sets $X_{n,\sigma}^- = \{x^- \in X_{(-\infty,0]} : \sigma \in \omega_1^+(x_{[-n.0]})\}$, that are open in $X_{(-\infty,0]}$. By property $g$ these open sets cover $X_{(-\infty,0]}$. There is a finite subcover $\{X_{n(\sigma):\sigma}^-, \sigma \in \Sigma\}$. With $n = \max\{n(\sigma) : \sigma \in \Sigma\}$ one has that the n-block system of $X$ is right-instantaneous. □

For a subshift $X \subset \Sigma^{\mathbb{Z}}$, for $L \in \mathbb{Z}_+$, and for mappings

$$\Psi^{(r)} : X_{[-L,L]} \to X_{[1,L+1]}$$

one formulates a condition

$$(RIa) : \Psi^{(r)}(a) \in \Gamma_{L+1}^+(x^-, a_{[-L,0]}), \qquad a \in X_{[-L,L]}, x^- \in \Gamma_\infty^-(a).$$

If a mapping $\Psi^{(r)} : X_{[-L,L]} \to X_{[1,L+1]}$ satisfies condition $(RIa)$ then for $0 \leq n < L$, and for $b^{(r)} \in X_{[-L-n,L]}$, the words $a_{n,b^{(r)}}^{(r)}$ that are given by

$$a_{n,b^{(r)}}^{(r)} = (b_{[-L-n,0]}^{(r)}, \Psi^{(r)}(b_{[-L,L]}^{(r)})_{(0,L-n]})$$

are in $X_{[-L-n,L-n]}$, and it is meaningful to impose on $\Psi^{(r)}$ a further condition

$$(RIb) : \Psi^{(r)}(a_{n,b^{(r)}}^{(r)}) = \Psi^{(r)}(b_{[-L-n,L-n]}^{(r)}), b^{(r)} \in X_{[-L-n,L]}, 0 \leq n < L.$$

We say that a mapping $\Psi^{(r)} : X_{[-L,L]} \to X_{[1,L+1]}$ that satisfies condition $(RIa)$ and also satisfies condition $(RIb)$ is an $RI$-mapping, and we say that a subshift that has an $RI$-mapping has property $RI$.

**Proposition 4.2.** *A subshift admits a presentation that has property g if and only if it has property IR.*

*Proof.* A subshift $X \subset \Sigma^{\mathbb{Z}}$ admits a right-instantaneous presentation if and only if it has property IR [8]. Apply Proposition 4.1. □

**Proposition 4.3.** *Let $X \subset \Sigma^{\mathbb{Z}}$ and $\widetilde{X} \subset \widetilde{\Sigma}^{\mathbb{Z}}$ be topologically conjugate subshifts, and let the subshift $\widetilde{X} \subset \widetilde{\Sigma}^{\mathbb{Z}}$ have property $(D)$. Then the subshift $X \subset \Sigma^{\mathbb{Z}}$ also has property $(D)$.*

*Proof.* A subshift has property $(D)$ if and only if one of its n-block systems has property $(D)$. To prove the proposition it is therefore sufficient to consider the situation that there is given a topological conjugacy $\varphi : X \to \widetilde{X}$ that is implemented by a 1-block map $\Phi : \Sigma \to \widetilde{\Sigma}$, with $\varphi^{-1}$ implemented for some $L \in \mathbb{Z}_+$, by a block map $\widetilde{\Phi}$ with coding window $[-L,L]$. Let there be given $a\sigma \in X_{[-I,0]}, I \geq 2L$. One has to find a $b \in X_{[-I-J,-I)}$, such that

$$b \in \Gamma^-(a), \quad \sigma \in \omega_1^+(ba). \tag{13}$$

For this, let

$$\widetilde{c} \in \widetilde{X}_{[-L,L]}, \widetilde{a} \in \widetilde{X}_{[-I-L,L]} \cap \Gamma^-(\widetilde{c}),$$

be such that

$$a\sigma = \widetilde{\Phi}(\widetilde{a}\widetilde{c}). \tag{14}$$

By property $(D)$ of $\widetilde{X}$ and by Lemma 3.3 there exists a

$$\widetilde{b} \in \widetilde{X}_{[-I-J-2L,-I-L)} \cap \Gamma^-(\widetilde{a}),$$

such that

$$\widetilde{c} \in \omega^+(\widetilde{b}_{[-I-J-L,-I)}, \widetilde{a}). \tag{15}$$

Then set

$$b = \widetilde{\Phi}(\widetilde{b}, \widetilde{a}_{[-I-L,-I+L)}) \tag{16}$$

and have by (14), (15) and (16) that (13) holds. □

An alternate proof of the invariance under topological conjugacy of property $(D)$ can be based on Nasu's theorem [13, Theorem 2.4] and on a notion of strong shift equivalence for measure $\lambda$-graph systems that is patterned after the notion of strong shift equivalence for $\lambda$-graph systems [12].

We describe prototype examples of subshifts with property $(D)$ and their presentations with property g. For this, we consider the Dyck inverse monoid with unit **1** and generating set $\{\alpha_\lambda, \alpha_\rho, \beta_\lambda, \beta_\rho,\}$, with relations

$$\alpha_\lambda \alpha_\rho = \beta_\lambda \beta_\rho = \mathbf{1}, \ \alpha_\lambda \beta_\rho = \beta_\lambda \alpha_\rho = 0.$$

The Dyck shift (on four symbols) is the subshift $D_2 \subset \{\alpha_\lambda, \alpha_\rho, \beta_\lambda, \beta_\rho,\}^\mathbb{Z}$ that contains all $x \in \{\alpha_\lambda, \alpha_\rho, \beta_\lambda, \beta_\rho,\}^\mathbb{Z}$ such that

$$\prod_{I_- \leq i < I_+} x_i \neq 0, \qquad I_-, I_+ \in \mathbb{Z}, I_- < I_+. \tag{17}$$

The Motzkin shift (on five symbols) is the subshift $M_2 \subset \{\mathbf{1}, \alpha_\lambda, \alpha_\rho, \beta_\lambda, \beta_\rho,\}^\mathbb{Z}$ that contains all $x \in \{\mathbf{1}, \alpha_\lambda, \alpha_\rho, \beta_\lambda, \beta_\rho,\}^\mathbb{Z}$ such that (17) holds. These presentations of the Dyck and Motzkin shifts have property $g$, and the Dyck and Motzkin shifts have property $(D)$. The Dyck and Motzkin shifts are prototypes of a class of subshifts that was introduced in [3] by giving presentations that were called $\mathcal{S}$-presentations. $\mathcal{S}$-presentations have property $g$ and subshifts that admit an $\mathcal{S}$-presentation have property $(D)$.